\documentclass{amsart}
\usepackage{graphicx}
\usepackage{verbatim}
\usepackage{wrapfig}

\usepackage{amsmath, amssymb, amsthm}
\usepackage{amsfonts}
\newtheorem{theorem}{Theorem}
\newtheorem{lemma}{Lemma}

\numberwithin{equation}{section}
\newtheorem{proposition}{Proposition}

\newtheorem{question}{Question}

\begin{document}

\begin{title}
{Bicircular signed-graphic matroids}
\end{title}

\author{Vaidy Sivaraman}
\address{Department of Mathematical Sciences, 
Library North 2200, 
Binghamton University,
Binghamton, New York 13902-6000.}
\email{vaidy@math.binghamton.edu}

\keywords{Matroid, bicircular, signed-graphic, transversal, series-parallel}
\date{June 26, 2013 \\  \text{      }  \text{  } \small {2010 Mathematics Subject Classification: 05C07}}   
\maketitle
\maketitle

\begin{abstract}

Several matroids can be defined on the edge set of a graph. Although historically the cycle matroid has been the most studied, in recent times, the bicircular matroid has cropped up in several places. A theorem of Matthews from late 1970s gives a characterization of graphs whose bicircular matroids are graphic. We give a characterization of graphs whose bicircular matroids are signed-graphic. \\

\end{abstract}

\begin{section}
{Introduction}

Matroids are combinatorial objects that were discovered when Whitney noticed strange similarities in properties between spanning trees of a connected graph and maximal linearly independent sets of columns of a matrix. Since then several classes of matroids have been studied in detail. One of the most important classes of matroids involves graphs, and consists of those matroids that can be realized as cycle matroids of graphs (defined later). Several results are known about the class of graphic matroids. In the last few decades a generalization of graphic matroids has proved to be an extremely important notion: the class of matroids arising from group-labelled graphs. The matroids that arise when the associated group has order two are called signed-graphic. In this paper we determine the intersection of two classes of matroids, namely, bicircular and signed-graphic. In late 1970s Matthews characterized graphs whose bicircular matroids are graphic. The work presented here is inspired from Matthews' paper \cite{
LRM} and his 
result can be obtained from the result here without any difficulty.

\end{section}

\begin{section}
{Terminology and notation}
Standard reference for matroid theory is Oxley's treatise \cite{JO}. $U_{k,n}$ denotes the uniform matroid with $n$ elements and rank $k$. $K_n$ denotes the complete graph on $n$ vertices. A subdivision of a graph $G$ is obtained by replacing some edges of $G$ by paths, where the internal vertices of the paths are disjoint from the vertices of $G$. Alternatively, a subdivision of $G$ is a graph obtained by repeatedly applying the operation of inserting a vertex of  degree two in  an edge. Two graphs $G$ and $H$ are said to be homeomorphic if there exists a graph $K$ such that both $G$ and $H$ are isomorphic to subdivisions of $K$. A matroid $N$ is said to be a minor of a matroid $M$ if a matroid isomorphic to $N$ can be obtained from $M$ by a sequence of deletions and contractions. A matroid is said to be binary if it is the vector matroid of a matrix with entries in $GF(2)$. A matroid is said to be ternary if it is the vector matroid of a matrix with entries in $ GF(3)$. 
The cycle matroid of a graph $G$, denoted $M(G)$, is the matroid on $E(G)$, where a set of edges is independent if the subgraph spanned by it contains no cycles. A matroid is graphic if it isomorphic to the cycle matroid of some graph. We use the symbol $\cong$ to denote isomorphism. \\
 
A signed graph $\Sigma$ is a pair $(G, \sigma)$ where $G$ is a graph and $\sigma : E(G) \rightarrow \{-1,1\}$ is a function. A circle (nonempty connected $2$-regular subgraph of $G$) is said to be positive if the product of the signs on its edges is $1$, and negative otherwise. The frame matroid of a signed graph $\Sigma$, denoted $M(\Sigma)$, is the matroid on the edge set of $\Sigma$, where a set of edges in $G$ is independent if the signed graph spanned by it has no positive circle and at most one negative circle in each component. A matroid is signed-graphic if it isomorphic to the frame matroid of some signed graph.\\
\end{section}

\begin{section}
{Statement of Matthews' theorem}
The bicircular matroid of a graph $G$, denoted $B(G)$, is the matroid on the edge set of $G$, where a set of edges in $G$ is independent if the subgraph spanned by it has at most one cycle in each component.

\begin{theorem}[Matthews, \cite{LRM}]
Let $G$ be a graph. Then the following conditions are equivalent:

\begin{enumerate}
\item $B(G)$ is graphic.
\item $B(G)$ is binary.
\item $B(G)$ is regular (that is, representable over every field).
\item Each component of $G$ can be obtained, by (repeated) addition of pendant edges, from either a theta graph or a graph homeomorphic to a tree with loops at some vertices.
\item $G$ has no subgraph homeomorphic to any of the graphs shown in  Figure 1  (where a graph with a dotted edge represents either the graph itself or the graph obtained when the dotted edge is contracted).
\end{enumerate}
\end{theorem}

\begin{figure}
\includegraphics[width=100mm, height=20mm]{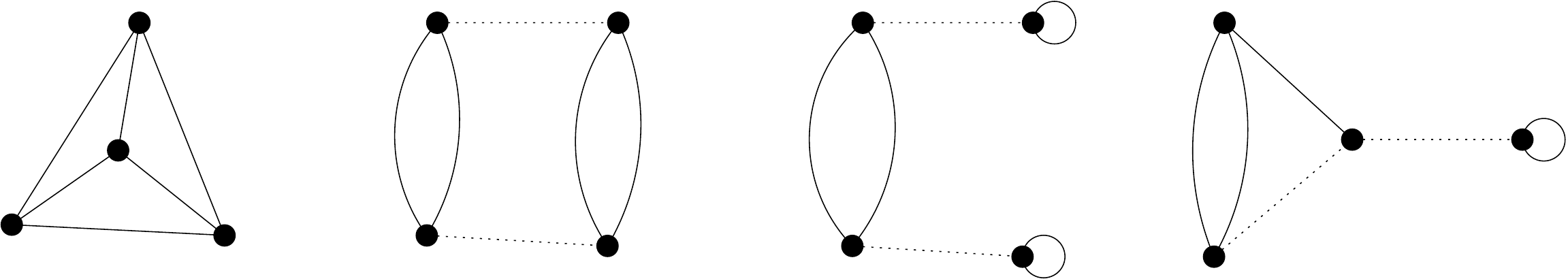}
\caption{Forbidden graphs in Matthews' characterization of graphs whose bicircular matroids are graphic.}
 \label{MatthewsGraphs}
\end{figure}

\end{section}

\begin{section}
{Statement of the main theorem}

  Let $k$ be a positive integer. $k$-skein is the graph with two vertices and $k$ edges, none of which is a loop.  A $k$-theta graph is a subdivision of  the $k$-skein. We are now ready to state the signed-graphic version of Matthews' theorem. 

\begin{theorem}\label{MAINRESULT}
Let $G$ be a graph. Then the following conditions are equivalent:

\begin{enumerate}
\item $B(G)$ is signed-graphic.
\item $B(G)$ is ternary.
\item $B(G)$ is near-regular (that is, representable over every field with at least three elements).
\item Each component of $G$ can be obtained, by (repeated) addition of pendant edges, from a subdivision of a tree with loops at some vertices and some edges doubled (an edge can be tripled (quadrupled) if one (both) of its endpoints is (are) pendant and loopless).
\item $G$ does not contain a subgraph that is a subdivision of any of  the graphs shown in Figure 2  (where a graph with a dotted edge represents either the graph itself or the graph obtained when the dotted edge is contracted).

\end{enumerate}
\end{theorem}

\begin{figure}[h]
\includegraphics[width=100mm, height=20mm]{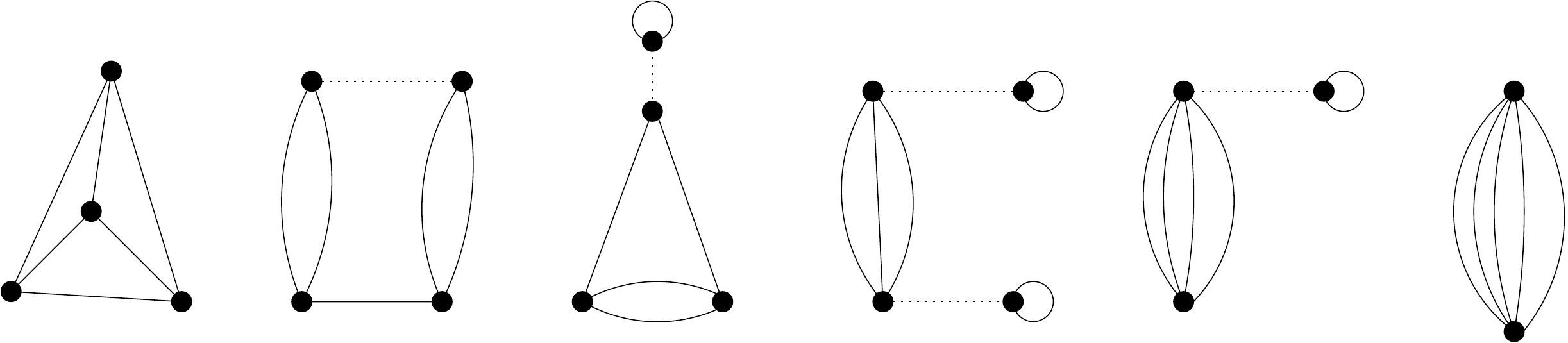}
\caption{Forbidden graphs in the characterization of graphs whose bicircular matroids are signed-graphic. We denote the graphs (from left to right) by $G_1,G_2,G_3,G_4,G_5$ and $G_6$.}
 \label{VaidyGraphs}
\end{figure}

\end{section}

\begin{section}
{Bicircular matroids and their properties}

Bicircular matroids were introduced by Sim\~{o}es-Pereira \cite{SP}, and studied in detail by Matthews \cite{LRM}.  \\

 The circuits of the bicircular matroid of a graph are edge sets of subgraphs that are subdivisions of the following three graphs: two vertices joined by three edges; two loops at the same vertex; two loops at distinct vertices that are joined by a single edge (called theta graphs, tight handcuffs, and loose handcuffs, respectively), because these are the only connected subgraphs which contain more than one cycle, but the deletion of any edge destroys all but one cycle. \\

 Bicircular matroids are transversal (Matthews \cite{LRM}), and hence inherit all nice properties of transversal matroids. In particular, they are representable over all sufficiently large fields and over all infinite fields (Piff and Welsh \cite{PW}). They are base-orderable, and can be written as a union of rank-1 matroids.  But their behavior is much more streamlined than that of transversal matroids: while the former is minor-closed, the latter is not. In fact, a recent project of DeVos et al. \cite{DGMR} is aimed at getting a complete list of excluded minors for the  class of bicircular matroids. \\

 It is important to note that, by definition, bicircular matroids are loopless, and hence it must be understood that when we take a minor of a bicircular matroid we do not allow contractions that produce loops. 

\end{section}

\begin{section}
{Signed-graphic matroids and their properties}

Signed-graphic matroids were introduced by Zaslavsky \cite{TZ1}, and studied in detail by Slilaty and Zaslavsky. A list of work done on topics related to signed graphs has been meticulously collected and maintained by Zaslavsky \cite{TZ2}.\\

The circuits of the frame matroid of a signed graph are edge sets of three types of subgraphs: positive circle, tight handcuff with both circles negative, and loose handcuff with both circles negative. For several other cryptomorphic definitions see \cite{TZ1}. \\

 The class of signed-graphic matroids properly contains the class of graphic matroids. Signed-graphic matroids need not be binary because the four-point line ($U_{2,4}$) is signed-graphic. Indeed, the signed graph consisting of a negative digon and two negative loops, one at each endpoint, is the unique signed graph whose frame matroid is $U_{2,4}$.  If an element of a signed-graphic matroid is deleted or contracted, the resulting matroid is still signed-graphic, and hence the class of signed-graphic matroids is closed under taking minors. However, the class of signed-graphic matroids is not closed under duality; for example, the bond matroid of the complete graph on $7$ vertices is not signed-graphic.  \\

 Signed-graphic matroids are dyadic, i.e,  representable over all fields whose characteristic is not  2 (Dowling-Zaslavsky \cite{TZ1}). The complete list of excluded minors for the  class of signed-graphic matroids is not known, although all such regular matroids are known (Qin, Slilaty, Zhou \cite{HSQ}). 
The intersection of the class of signed-graphic matroids with the class of cographic matroids has been characterized by Slilaty \cite{DS}.

\end{section}
\begin{section}
{Proof of the main theorem}

\begin{proof}[Proof of Theorem \ref{MAINRESULT}]
We will prove the following chain of five implications: \\

 $ (1) \Rightarrow (2) \Rightarrow (3) \Rightarrow  (5) \Rightarrow (4) \Rightarrow (1).$ \\

\begin{itemize}

\item $(1) \Rightarrow (2)$ It is well known that signed-graphic matroids are ternary (Dowling-Zaslavsky, \cite{TZ1}). A suitably chosen incidence matrix is a ternary representation for the frame matroid. \\

\item  $(2) \Rightarrow (3)$ Bicircular matroids are representable over all sufficiently large finite  fields and over all infinite fields (in fact, this is true for the bigger class of transversal matroids) (Piff and Welsh, \cite{PW}). In particular, a bicircular matroid is representable over $\mathbb{Q}$ and a field of characteristic $2$. A theorem of  Whittle \cite{GW} says that if a matroid $M$ is representable over $GF(3), \mathbb{Q}$, and over a field of characteristic two, then $M$ is near-regular. Suppose $B(G)$ is ternary. We know that $B(G)$ is representable over $\mathbb{Q}$ and over a field of characteristic $2$. Invoking Whittle's theorem, we see that $B(G)$ is near-regular.   \\

\item $(3) \Rightarrow (5)$
The following statements are easy to check:
\begin{itemize}
\item $B(G_1) \cong U_{4,6}$. 
\item $U_{3,5}$ is a minor of $B(G_2)$ and $B(G_3)$. 
\item $U_{2,5}$ is a minor of $B(G_4)$ and $B(G_5)$. 
\item $B(G_6) \cong U_{2,5}$. 
\end{itemize}

 Suppose $B(G)$ is near-regular. In particular,  $B(G)$ is ternary.  But none of $U_{4,6}, U_{2,5},$ and $U_{3,5}$ is ternary; in fact $U_{2,5}$ and $U_{3,5}$ are excluded minors for the class of ternary matroids.  This together with the fact that the class of bicircular matroids is closed under taking minors shows that $G$ has no subgraph homeomorphic to any of the graphs $G_i$. \\

\item $(5) \Rightarrow (4)$  

 We need the following well-known characterization of graphs none of whose subgraphs is a subdivision of $K_4$. 

\begin{lemma}\label{SPLEMMA}
Let $G$ be a connected graph. The following statements are equivalent:

\begin{itemize}
\item No subgraph of $G$ is a subdivision of $K_4$.
\item Every component of $G$ can be iteratively constructed from $K_1$ by the following operations:
  \begin{itemize}
  \item Adding a loop.
  \item Adding a pendant edge.
  \item Adding an edge in parallel to an existing edge.
  \item Subdividing an edge.
  \end{itemize}
\end{itemize}
\end{lemma}

A graph none of whose subgraphs is a subdivision of $K_4$ is called a series-parallel graph. A proof of the above lemma can be easily obtained from the following easy result, whose proof can be found in many standard graph theory texts:

\begin{proposition}
Let $G$ be a (non-null) simple graph with minimum degree at least 3. Then $G$ contains a subgraph that is a subdivision of $K_4$. 
\end{proposition}

By virtue of Lemma \ref{SPLEMMA}, we can restrict attention to graphs, all of whose components can be obtained from $K_1$ by a sequence of the four operations mentioned above.

\begin{lemma}\label{THREETHETA}
Let $G$ be obtained by adding an edge $f$ in parallel to an edge $e$ in a $3$-theta graph $H$. Then $G$ is either a subdivision of $G_2$ or a $4$-theta graph. 
\end{lemma}

\begin{proof}
Let $u,v$ be the trivalent vertices of $H$. If the endpoints of $e$ are different from $u,v$, then $G$ is a subdivision of $G_2$ where the dotted edge is an edge. If $f$ has exactly one endpoint in $u,v$, then $G$ is a subdivision of $G_2$ where the dotted edge is contracted. If $u,v$ are the endpoints of $f$, then $G$ is a $4$-theta graph. \end{proof}

\begin{lemma}\label{TWOSKEIN}
Let $G$ be a series-parallel graph with minimum degree at least $3$. Suppose that no subgraph of $G$ is a subdivision of either $G_2$ or $G_3$. Let $B$ be a block in which the base edge is duplicated exactly once. Then $B$ is isomorphic to the $2$-skein. 
\end{lemma}

\begin{proof}
Let the endpoints of the base edge be $u$ and $v$. Suppose that the base edge has an internal vertex $w$. Up to symmetry, there are three possibilities: 
(a) Two edges with endpoints $u$ and $w$, and two edges with endpoints $v$ and $w$. Together with the duplicated edge $f$, we have $G_2$.   
(b) There exists a vertex $x$ and two edges with endpoints $u$ and $w$, and two edges with endpoints $v$ and $x$ and an edge with endpoints $w$ and $x$. Together with the duplicated edge $f$, we have $G_2$.  
(c) There is another block containing $v$. Hence there is an endblock containing a cycle and we have a subdivision of $G_3$. 
We conclude that the base edge has no internal vertex and hence $B$ is isomorphic to the $2$-skein. 
\end{proof}

\begin{lemma}\label{SKEINLEMMA}
Let $G$ be a series-parallel graph with minimum degree at least $3$. Suppose that no subgraph of $G$ is a subdivision of either $G_2$ or $G_3$. Then every block of $G$ is either a loop or a $k$-skein for some positive integer $k$. 
\end{lemma}

\begin{proof}
Let $B$ be a block of $G$. We may assume that $B$ is not a loop. By Lemma \ref{THREETHETA}, we know that if the base edge is duplicated at least twice, it has to be a theta-graph. If the base edge is duplicated only once, by Lemma \ref{TWOSKEIN}, $B$ is isomorphic to the $2$-skein. If the base edge is not duplicated, then $B$ just consists of the single edge, and is isomorphic to a $1$-skein. 
\end{proof}

The following lemma can be easily established using similar arguments. 

\begin{lemma}\label{ENDBLOCKLEMMA}
Let $G$ be a series-parallel graph with minimum degree at least $3$. Suppose that no subgraph of $G$ is a subdivision of $G_4$. If $B$ is a block of $G$ isomorphic to the $3$-skein, then $B$ is an endblock. 
\end{lemma}

Armed with the above lemmas, we are now ready to prove the implication $(5) \Rightarrow (4)$ . Suppose no subgraph of $G$ is a subdivision of any of $G_i$, $1 \leq i \leq 6$. Let $G'$ be obtained from $G$ by a sequence of operations, the operations being deleting pendant vertices and contracting edges in series. It is routine to check that $G'$ is well-defined. \\

It is important to note that no subgraph of $G'$ is a subdivision of any of $G_i$, $1 \leq i \leq 6$. In particular, $G'$ is series-parallel and none of its subgraph is a subdivision of $G_2$ or $G_3$. By Lemma \ref{SKEINLEMMA}, each block of $G'$ is either a loop or a $k$-skein. Note that $G_6$ is the $5$-skein, hence each non-loop block of $G'$ is a $k$-skein for $k = 1,2,3,$ or $4$. Suppose $G'$ contains a block that is isomorphic to the $4$-skein. Since $G'$ does not contain a subdivision of $G_5$, we conclude that $G'$ is isomorphic to the $4$-skein. Now, suppose that no block of $G'$ is isomorphic to the $4$-skein. Let $G'$ have a block $B$ that is isomorphic to the $3$-skein. Then, by Lemma \ref{ENDBLOCKLEMMA}, $B$ is an endblock of $G'$.  Hence $G'$ is isomorphic either to the $4$-skein or to a tree with some edges doubled and loops added, and tripling of edges allowed only when one of its endpoints is a leaf of the tree and has no loops on it. By reversing the operations done to get $G'$ from 
$G$, we see that $G$ can be obtained, by (repeated) addition of pendant edges, from a subdivision of a tree with loops at some vertices and some edges doubled (an edge can be tripled (quadrupled) if one (both) of its endpoints is (are) pendant and loopless); \\

\item $(4) \Rightarrow (1)$: 
Since the class of signed-graphic matroids is closed under taking direct sums, it suffices to prove the result for connected graphs. Henceforth we will assume that $G$ is connected. Note that the bicircular matroid of the 3-skein is isomorphic to the frame matroid of a positive triangle, and the bicircular matroid of the 4-skein or the 3-skein with a loop is isomorphic to the four-point line ($U_{2,4}$) which is isomorphic to the signed graph consisting of a negative digon and two negative loops, one at each endpoint.   

Suppose $G$ is a graph obtained from a tree by doubling some edges and adding loops. We will construct a signed graph $\Sigma$ such that $B(G) = M(\Sigma)$, proving that $B(G)$ is signed-graphic. Let the underlying graph of $\Sigma$ be $G$. Also, make all loops and exactly one edge of each digon negative. Note that, because of the special structure of $G$, a cycle in $G$ corresponds to a negative circle in $\Sigma$. Also, there are no positive circles in $\Sigma$. Hence  $B(G) = M(\Sigma)$. \\

Suppose $G$ is a graph obtained from a tree by doubling some edges and adding loops, and some pendant edges tripled when the corresponding pendant vertex has no loops. Then we can construct a graph $G'$ by replacing an edge in each triple (three edges in parallel) by a loop at the corresponding pendant endpoint. It is easy to check that $B(G) = B(G')$. Now we repeat the procedure above to get a signed graph $\Sigma$ with $B(G') = M(\Sigma)$. This implies $B(G) = M(\Sigma)$ and  hence $B(G)$ is signed-graphic. \\

We conclude, by using the following two easy lemmas (whose trivial proofs we omit) to show, that if $G$ is of the form mentioned in (4), then $B(G)$ is signed-graphic. \\

\begin{lemma}
Let $H$ be obtained from a graph $G$ by subdividing an edge. If $B(G)$ is signed-graphic, so is $B(H)$.
\end{lemma}

\begin{lemma}
Let $H$ be obtained from a graph $G$ by adding a pendant edge. If $B(G)$ is signed-graphic, so is $B(H)$.
\end{lemma}

\end{itemize}
This concludes all the five implications, and hence concludes the proof of the main theorem.
\end{proof}

\end{section}

\begin{section}
{Questions and concluding remarks}

We conclude with some open-ended questions that are similar in spirit to the content of this paper.

\begin{question}
Characterize signed graphs whose matroids are bicircular.
\end{question}

\begin{question}
Characterize signed graphs whose matroids are transversal.
\end{question}

\begin{question}
Characterize set systems whose transversal matroids are signed-graphic.
\end{question}

\begin{question}
When is the dual of a signed-graphic matroid signed-graphic? 
\end{question}

\begin{question}[Welsh \cite{DJAW}]
When is the dual of a transversal matroid transversal? 
\end{question}

\begin{question}
It is known that we cannot determine in polynomial time whether a matroid, given in terms of independence oracle, is signed-graphic (Geelen-Mayhew \cite{MWZ}).  It is also known that we cannot determine in polynomial time whether a matroid, given in terms of independence oracle, is bicircular. Can we determine in polynomial time whether a matroid, given in terms of independence oracle, is signed-graphic and bicircular? 
\end{question}

A graph can be thought of as a signed graph with all edges positive. Since a graph with all edges positive has no negative circles, the frame matroid of an all-positive graph is equal to the cycle matroid of the graph, and hence graphic matroids are signed-graphic. Using the fact that $U_{2,4}$ is not graphic, Matthews' theorem can be obtained as an immediate corollary of the main result of this paper. 

\end{section}

\begin{section}
{Acknowledgements}
The work presented here is part of my PhD dissertation \cite{VS} written under the guidance of Neil Robertson at The Ohio State University. I would like to thank James Oxely, Neil Robertson, Daniel Slilaty and Thomas Zaslavsky for their help in preparing this article. Luis Goddyn gave an excellent talk on the excluded minors for the class of bicircular matroids at the Third Workshop on Graphs and Matroids at Maastricht (organized by Bert Gerards), and that was the starting point for this work. 
\end{section}

\end{document}